\documentclass[fleqn]{mat01}
\usepackage{times,mathtimy,amssymb,latexsym,amsbsy}
\begin{document}

\newtheorem{lem}{\it Lemma}
\newtheorem{theor}{\bf Theorem}

\def\conremar{\trivlist\item[\hskip\labelsep{{\it Concluding Remarks.}}]}

\title{Zeta function of the projective curve $\pmb{aY^{2\,l} = bX^{2\,l} +
cZ^{2\,l}}$ over a class of finite fields, for odd primes $\pmb{l}$}

\markboth{N Anuradha}{Zeta function of the projective curve $aY^{2l} =
bX^{2l} + cZ^{2l}$}

\author{N ANURADHA}

\address{Institute of Mathematical Sciences, C.I.T. Campus, Taramani,
Chennai~600~113, India\\
\noindent E-mail: anuradha@imsc.res.in}

\volume{115}

\mon{February}

\parts{1}

\pubyear{2005}

\Date{MS received 18 June 2003}

\begin{abstract}
Let $p$ and $l$ be rational primes such that $l$ is odd and the order of
$p$ modulo $l$ is even. For such primes $p$ and $l$, and for $e=l, 2l$,
we consider the non-singular projective curves $aY^e = bX^e + cZ^e$
($abc \neq 0$) defined over finite fields $\mathbf{F}_q$ such that
$q=p^\alpha \equiv 1(\bmod~{e})$. We see that the Fermat curves correspond
precisely to those curves among each class (for $e=l,2l$), that are
maximal or minimal over $\mathbf{F}_q$. We observe that each Fermat
prime gives rise to explicit maximal and minimal curves over finite
fields of characteristic 2. For $e=2l$, we explicitly determine the
$\zeta$-function(s) for this class of curves, over $\mathbf{F}_q$, as
rational functions in the variable $t$, for distinct cases of $a,b$, and
$c$, in $\mathbf{F}_q^*$. The $\zeta$-function in each case is seen to
satisfy the Weil conjectures (now theorems) for this concrete class of
curves.

For $e=l, 2l$, we determine the class numbers for the function fields
associated to each class of curves over $\mathbf{F}_q$. As a
consequence, when the field of definition of the curve(s) is fixed, this
provides concrete information on the growth of class numbers for
constant field extensions of the function field(s) of the curve(s).
\end{abstract}

\keyword{Finite fields; curves; maximal curves; zeta functions; function
fields.}

\maketitle

\section{Introduction}\label{sec1}

Let $p$ and $l$ be rational primes such that $l$ is odd and the order of
$p$ modulo $l$, written as ${\rm ord}\ p(\bmod~{l})$, is even. Let $f =
{\rm ord}\ p(\bmod~{l})$; this is defined to be the least positive
integer such that $p^f \equiv 1(\bmod~{l})$. For such primes $p$ and
$l$, we consider finite fields $\mathbf{F}_q$ such that $q = p^\alpha
\equiv 1(\bmod~{e})$, for $e=l,2l$; thus $\alpha = fs$ for some integer
$s \geq 1$. If $e=2l$, clearly $p$ is odd.

In (\cite{an-1}, Theorems~6, 7), we had considered the non-singular
projective curves $aY^e = bX^e + cZ^e$ ($abc \neq 0$, and $e=l, 2l$)
defined over such finite fields $\mathbf{F}_q$, and had explicitly
obtained the number of $\mathbf{F}_{q^n}$-rational points on these
curves for each integer $n \geq 1$. This was done by applying explicit
results obtained in \cite{an-1} for the cyclotomic numbers of order $e$
over $\mathbf{F}_q$.

Further, for the case $e=l$, we had obtained in (\cite{an-1}, Theorem~8)
the explicit $\zeta$-function(s) for this class of curves defined over
$\mathbf{F}_q$. In this paper, we consider the case $e=2l$ and apply the
results of (\cite{an-1}, Theorem~7) to obtain the explicit
$\zeta$-function(s), in Theorem~\ref{thm-1} (\S\ref{sec3}), for the
class of non-singular projective curves $aY^{2l} = bX^{2l} + cZ^{2l}$
$(abc \neq 0)$ defined over $\mathbf{F}_q$, as rational function(s) in
the variable $t$. We do this for all distinct cases of $a,b$, and $c$,
in $\mathbf{F}_q^*$. There are seven distinct cases, and the
$\zeta$-function in each case is seen to satisfy the Weil conjectures
(proven in generality) for this concrete class of curves.

In \S\ref{sec2}, we define maximal and minimal curves over finite
fields, and we interpret the results obtained in (\cite{an-1},
Theorems~6, 7), in this context, to link these results with facts
previously known in the literature. In addition, we make some simple but
pertinent observations pertaining to these results.

In \S\ref{sec3}, as a consequence to the explicit $\zeta$-functions
obtained in Theorem~\ref{thm-1} and (\cite{an-1}, Theorem 8), for the
projective curves $aY^e = bX^e + cZ^e$ ($e=l, 2l$) defined over
$\mathbf{F}_q$, we obtain the class numbers of the associated function
fields in Theorems~\ref{thm-2} and \ref{thm-3}, for all distinct cases
of $a,b,c \in \mathbf{F}_q^*$. Further, for $e=l,2l$, if we fix the
field of definition $\mathbf{F}_q$, and consider the curve(s) over all
finite extensions of $\mathbf{F}_q$, these results provide concrete
information on the growth of class numbers for constant field extensions
of the function field of the curve(s) over $\mathbf{F}_q$.

For easy reference, we restate (\cite{an-1}, Theorems~6, 7) as
Lemmas~\ref{lem-1} and~\ref{lem-2} below:

\begin{lem}\label{lem-1}\hskip -.2pc${\rm ([1], Theorem~6).}$\ \ Let $p$ be
any prime such that $f = {\rm ord}\ p(\bmod~{l})$ is even. Let $q =
p^\alpha \equiv 1(\bmod~{l})${\rm ,} and $\alpha = fs$ for some integer $s
\geq 1$. Consider the projective curve $aY^l = bX^l + cZ^l$ $(abc \neq
0)$ defined over the finite field $\mathbf{F}_q$. Fix any generator
$\gamma$ of $\mathbf{F}_q^*$ and let ${\rm ind}_\gamma(b/c) \equiv
i(\bmod~{l})$ and ${\rm ind}_\gamma(a/c) \equiv j(\bmod~{l})$. Then for
each $n \geq 1${\rm ,} the number $a_l(n)$ of $\mathbf{F}_{q^n}$-rational
points on this curve is given as below{\rm :}
\begin{align*}
a_l(n) &= q^n + 1 -\! (l-1)(l-2)(-1)^{ns}q^{n/2},\quad if~
~in, jn \equiv 0(\bmod~{l}), \\[.2pc]
a_l(n) &= q^n + 1 - 2(-1)^{ns}q^{n/2},\quad if~
~in,jn,in-\!jn \not\equiv 0(\bmod~{l}),  \\[.2pc]
a_l(n) &= q^n + 1 + (l-2)(-1)^{ns}q^{n/2},\quad
in~all~other~cases~of \\
&\hskip 4.9cm in,~jn(\bmod~{l}).
\end{align*}
\end{lem}

\begin{lem}\label{lem-2}\hskip -.2pc${\rm ([1], Theorem~7).}$\ \
Let $p$ be an odd prime such that $f = {\rm ord}\ p(\bmod~{l})$ is even.
Let $q = p^\alpha \equiv 1(\bmod~{2l})${\rm ,} and $\alpha = fs$ for some
integer $s \geq 1$. Consider the projective curve $aY^{2l} = bX^{2l} +
cZ^{2l}$ $(abc \neq 0)$ defined over the finite field $\mathbf{F}_q$.
Fix any generator $\gamma$ of $\mathbf{F}_q^*$ and let ${\rm
ind}_\gamma(b/c) \equiv i(\bmod~{2l})$ and ${\rm ind}_\gamma(a/c) \equiv
j(\bmod~{2l})$. Then for each $n \geq 1${\rm ,} the number $a_{2l}(n)$ of
$\mathbf{F}_{q^n}$-rational points on this curve is given as\break below{\rm
:}
\begin{align*}
a_{2l}(n) &= q^n + 1 - (2l-1)(2l-2)(-1)^{ns}q^{n/2},\quad
if~~in,jn \equiv 0(\bmod~{2l}),\\[.2pc]
a_{2l}(n) &= q^n + 1 -2(-1)^{ns}q^{n/2},\quad if~~in,jn,in-jn
\not\equiv 0(\bmod~{2l}),\\[.2pc]
a_{2l}(n) &= q^n + 1 + 2(l-1)(-1)^{ns}q^{n/2},\quad
in~all~other~cases~of \\
&\hskip 5.1cm in,~jn(\bmod~{2l}).
\end{align*}
\end{lem}

\section{Maximal curves defined over finite fields}\label{sec2}

Denote by $C/\mathbf{F}_q$ a non-singular projective algebraic curve $C$
defined over a finite field $\mathbf{F}_q$. Denote by $a(n,C)$ the
number of $\mathbf{F}_{q^n}$-rational points on $C$, for each integer $n
\geq 1$. Let $g \geq 0$ be the genus of $C$. The Weil conjectures for
the curve $C$ state that the $\zeta$-function of $C$ over $\mathbf{F}_q$, which
is defined as
\begin{equation*}
Z(t,C) = Z(t, C/\mathbf{F}_q) = \exp \left(\sum_{n=1}^{\infty}
\frac{a(n,C) t^n}{n} \right),
\end{equation*}
satisfies the following properties:
\begin{enumerate}
\renewcommand\labelenumi{\rm \arabic{enumi}.}
\leftskip -.2pc
\item $Z(t,C)$ is a rational function in the variable $t$ of the form
${P(t)}/{(1-t)(1-qt)}$, where $P(t)$ is a polynomial in $t$ of degree
$2g$, having integer coefficients, leading term $q^g$, and constant term
$1$.

\item $Z(t,C)$ satisfies a functional equation given by
\begin{equation*}
\hskip -1.25pc q^{g-1}t^{2g-2} Z(1/qt, C) = Z(t, C).
\end{equation*}
Equivalently, if we express $P(t) = \prod_{k=1}^{2g}(1-\alpha_k t)$, we
may pair the $\alpha_k$'s in such a way that $\alpha_k \alpha_{g+k} = q$
for $1 \leq k \leq g$.

\item The reciprocal roots $\alpha_k$ of $P(t)$ satisfy the property
that $|\alpha_k| = q^{1/2}$ for $1 \leq k \leq 2g$. This is known as the
Riemann hypothesis for $C/\mathbf{F}_q$.
\end{enumerate}

These conjectures were first stated (in full generality) by Andr\'e Weil
\cite{w-2} in 1949, for non-singular projective varieties of dimension
$\geq 1$ defined over finite fields. Curves are varieties of dimension
1. These conjectures have been proven in complete generality (see, for
example, \cite{freit}). For a proof of these conjectures for curves, see
\cite{cmor}.

The first general proof of these conjectures for {\it curves} was given
by Weil \cite{w-1}. He showed that for such a curve $C/\mathbf{F}_q$,
\begin{equation*}
a(n, C) = q^n + 1 - \sum_{k=1}^{2g} \alpha_k^n, \quad \text{for each } n
\geq 1.
\end{equation*}
As a consequence of the Riemann hypothesis for $C/\mathbf{F}_q$, he
obtained the following bounds on $a(n,C)$, given by
\begin{equation*}
|a(n,C) - (q^n + 1)| \leq 2g q^{n/2}, \quad \text{for each } n \geq 1.
\end{equation*}
These bounds, and a proof of the Riemann hypothesis, were earlier
obtained by Hasse in 1936 for curves of genus $g=1$ (or elliptic
curves), and have come to be known as the Hasse--Weil bounds for the
curve $C$, sometimes simply referred to as the Weil bounds\break for~$C$.

There has been considerable interest and search in the literature for
curves $C/\mathbf{F}_q$ for which the upper Weil bounds are attained for
the number of $\mathbf{F}_q$-rational points on $C$. Such curves are
called maximal, and the associated function fields are called maximal
function fields. Maximal curves are of theoretical interest; they
provide examples of curves with large automorphism groups, and have
interesting arithmetic and geometric properties (cf.
\cite{fgt,aik,leo,r-sti,sti-1}). Such curves, and curves
$C/\mathbf{F}_q$ with large number of $\mathbf{F}_q$-rational points,
find important applications in coding theory, since the construction by
Goppa of codes with good parameters from such curves (cf.
\cite{go-1,go-2}).

In keeping with the terminology for maximal curves, one may define
minimal curves to be curves $C/\mathbf{F}_q$ for which the lower Weil
bounds are attained for the number of $\mathbf{F}_q$-rational points
on $C$ (i.e, $a(1,C) = q + 1 - 2g q^{1/2}$).  It is clear from the
expression for the Weil bounds that curves $C/\mathbf{F}_q$ are
maximal or minimal only when $q$ is a square (even power of $p$),
or the genus $g=0$.

As a special case of the curves treated in Lemmas~\ref{lem-1} and
\ref{lem-2}, consider the Fermat curves $Y^e = X^e + Z^e$ (for $e=l,2l$)
defined over finite fields $\mathbf{F}_q, \: q = p^\alpha \equiv
1(\bmod~{e})$, when $f = {\rm ord}\ p(\bmod~{l})$ is even. If we fix $q
= p^f$, it is clear from Lemmas~\ref{lem-1} and~\ref{lem-2} that these
curves are maximal over finite odd degree extensions of $\mathbf{F}_q$,
and are minimal over finite even degree extensions of $\mathbf{F}_q$. It
would thus appear that there is a close inter-relationship between
maximal and minimal curves defined over finite fields.

Further, keeping $q=p^f$, if we write $q-1 = et$, for $e=l,2l$, and $t
\geq 1$, then since $f = {\rm ord}\ p(\bmod~{l})$, it follows that
$q^{1/2} + 1 = et'$ for some $t' \geq 1, t'| t$. For $t'=1$, the
corresponding Fermat curves are then defined by
\begin{equation*}
Y^{q^{1/2} + 1} = X^{q^{1/2} + 1} + Z^{q^{1/2} + 1}
\end{equation*}
over the finite field $\mathbf{F}_q$. These are just the Hermitian
curves which have been studied in the literature and known to be maximal
over $\mathbf{F}_q$. The corresponding function fields are called
Hermitian function fields. Hermitian curves have been characterized as
the (essentially) unique maximal curves over $\mathbf{F}_q$ with genus
$g = q^{1/2}(q^{1/2}-1)/2$. This is the maximum possible genus for a
maximal curve defined over $\mathbf{F}_q$ (cf. \cite{ihara,r-sti}). For
$t' > 1$, we have $e| q^{1/2} + 1$, and the corresponding Fermat curves
are again known to be maximal over $\mathbf{F}_q$; these are not
Hermitian, but the function fields associated to these curves occur as
subfields of the Hermitian function field (cf. (\cite{sti-2},
pp.~196--203)).

The case when $l$ is a Fermat prime is interesting; if $l = 2^{2^n} +
1$, the corresponding Fermat curve $Y^l = X^l + Z^l$ is a Hermitian
curve over the finite field $\mathbf{F}_q, \: q = 2^{2^{n+1}}$, of
characteristic 2, with genus $g = 2^{2^n -1}(2^{2^n} -1)$. Further, as
observed above, these curves are maximal (resp. minimal) over finite odd
degree (resp. even degree) extensions of $\mathbf{F}_q$. Thus each
Fermat prime gives rise to explicit maximal and minimal curves over
finite fields of characteristic 2. The converse, however, is not true.
For example, take $n =5$; then $r = 2^{2^5} + 1$ is {\it not} a prime,
but the corresponding Fermat curve $Y^r = X^r + Z^r$ is maximal (or
minimal) over finite extensions of the field $\mathbf{F}_q, \: q =
2^{2^6}$.

From the explicit results in Lemmas~\ref{lem-1} and~\ref{lem-2}, it is
also clear that the only class of coefficients $a,b,c$ in
$\mathbf{F}_q^*$ for which the curves $aY^e = bX^e + cZ^e$ are maximal
(or minimal) over $\mathbf{F}_q$ are those that correspond to the Fermat
curves $Y^e = X^e + Z^e$ (for $e=l, 2l$) (i.e., the coefficients reduce
to the case $a=b=c=1$). The cases for $a,b,c$ which do not correspond to
the Fermat curves are {\it never} maximal or minimal. (Note that
in Lemma~\ref{lem-1}, for $l=3$, we have $f = {\rm ord}\ p(\bmod~{3}) =
2$, and $q=p^{2s}$. For $s$ odd, each element of $\mathbf{F}_{p^s}^*$ is
a cube, and hence, all cases when $a,b,c \in \mathbf{F}_{p^s}^*$
correspond to the Fermat curve $Y^3 = X^3 + Z^3$, and this curve is
maximal over $\mathbf{F}_q$.)

\section{Zeta function(s) of the projective curve $\pmb{aY^{2\,l} = bX^{2\,l} +
cZ^{2\,l}}$ over $\mathbf{F}_{\pmb{q}}$}\label{sec3}

\begin{theor}[\!]\label{thm-1}
Let $p$ and $l$ be odd rational primes such that $f = {\rm ord}\
p(\bmod~{l})$ is even. Consider the projective curve $C \colon aY^{2l} =
bX^{2l} + cZ^{2l}$ $(abc \neq 0)$ defined over the finite field
$\mathbf{F}_q${\rm ,} where $q = p^\alpha \equiv 1(\bmod~{2l})${\rm ,}
and $\alpha = fs$ for $s \geq 1$. Fix a generator $\gamma$ of
$\mathbf{F}_q^*$ and let ${\rm ind}_{\gamma}(b/c) \equiv i(\bmod~{2l})$
and ${\rm ind}_{\gamma}(a/c) \equiv j(\bmod~{2l})$. Let $\theta =
(-1)^sq^{1/2}$ and let $\zeta$ be any primitive {\rm (}complex{\rm )} $l$-th root of
unity. Then the $\zeta$-function $Z(t,C)$ of the curve $C/\mathbf{F}_q$
is a rational function in the variable $t${\rm ,} of the form
$P(t)/(1-t)(1-qt)${\rm ,} and the polynomial $P(t)$ is given explicitly
for distinct cases of $i,j(\bmod~{2l})$ as below{\rm :}
\begin{enumerate}
\renewcommand\labelenumi{\rm \arabic{enumi}.}
\leftskip -.2pc
\item  For $i \equiv j \equiv 0(\bmod~{2l})${\rm ,}
\begin{equation*}
\hskip -1.25pc P(t) = (1- \theta t)^{(2l-1)(2l-2)}.
\end{equation*}

\item For $i,j \equiv 0(\bmod~{2}),\, i,j, i-j \not\equiv
0(\bmod~{2l})${\rm ,}
\begin{equation*}
\hskip -1.25pc P(t) = (1 - \theta t)^{4l-4} \prod_{r=1}^{l-1} (1 -
\zeta^r \theta t)^{4l-6}.
\end{equation*}

\item For $i, j \equiv 0(\bmod~{2})$ and {\rm (i)} $i \equiv 0,\, j
\not\equiv 0(\bmod~{2l})${\rm ,} {\rm (ii)} $i \not\equiv 0,\, j \equiv
0(\bmod~{2l})${\rm ,} {\rm (iii)}~$i,j \not\equiv 0,\, i \equiv
j(\bmod~{2l})${\rm ,}
\begin{equation*}
\hskip -1.25pc P(t) = (1 - \theta t)^{2l-2} \prod_{r=1}^{l-1}(1 -
\zeta^r \theta t)^{4l-4}.
\end{equation*}

\item For {\rm (i)} $i \equiv 0,\, j \equiv l(\bmod~{2l})${\rm ,} {\rm
(ii)} $i \equiv l,\, j \equiv 0(\bmod~{2l})${\rm ,} and {\rm (iii)} $i
\equiv j \equiv l(\bmod~{2l})${\rm ,}
\begin{equation*}
\hskip -1.25pc P(t) = (1 - \theta t)^{2(l-1)^2} (1 + \theta t)^{2l(l-1)}.
\end{equation*}

\item For {\rm (i)} $j \not\equiv 0(\bmod~{2}),\, i \equiv 0,\, j
\not\equiv l(\bmod~{2l})${\rm ,} {\rm (ii)} $i \not\equiv
0(\bmod~{2}),\, i \not\equiv l,\, j \equiv 0(\bmod~{2l})${\rm ,} and
{\rm (iii)} $i \not\equiv 0(\bmod~{2}),\, i \equiv j,\, i \not \equiv
l(\bmod~{2l})${\rm ,}
\begin{align*}
\hskip -1.25pc P(t) &= (1 + \theta t)^{2l-2} \prod_{r=1}^{l-1}((1 -
\zeta^r \theta t)(1 + \zeta^r \theta t))^{2l-2}\\[.3pc]
\hskip -1.25pc &= \prod_{r=1}^{2l-1}(1 - \xi^r \theta t)^{2l-2},
\end{align*}
where $\xi$ is a primitive complex $2l$-th root of unity.

\item For $i,j, i-j \not\equiv 0,l(\bmod~{2l})$ and {\rm (i)} $i
\not\equiv j(\bmod~{2})${\rm ,} {\rm (ii)} $i,j \not\equiv
0(\bmod~{2})${\rm ,}
\begin{equation*}
\hskip -1.25pc P(t) = ((1 - \theta t)(1 + \theta t))^{2l-2}
\prod_{r=1}^{l-1}((1 - \zeta^r \theta t)^{2l-4} (1 + \zeta^r \theta
t)^{2l-2}).
\end{equation*}

\item For {\rm (i)} $i \equiv l,\, j \not\equiv 0,l(\bmod~{2l})${\rm ,}
{\rm (ii)} $i \not\equiv 0,l,\, j \equiv l(\bmod~{2l})${\rm ,} and {\rm
(iii)} $i,j \not\equiv 0,l,\, i-j \equiv l(\bmod~{2l})${\rm ,}
\begin{equation*}
\hskip -1.25pc P(t) = ((1 - \theta t)(1 + \theta t))^{l-1}
\prod_{r=1}^{l-1}((1 - \zeta^r \theta t)^{2l-3} (1 + \zeta^r \theta
t)^{2l-1}).
\end{equation*}
\end{enumerate}
\end{theor}

\begin{proof} The number $a_{2l}(n)$ of $\mathbf{F}_{q^n}$-rational
points on the curve $C$, for each $n \geq 1$, has been determined
explicitly in (\cite{an-1}, Theorem 7) (cf. Lemma~\ref{lem-2} above).
Taking into account the distinct cases that arise when \hbox{$l| n, l
\!\not|\, n, 2 | n$,} and \hbox{$2 \!\not|\, n$,} and substituting the
corresponding values for $a_{2l}(n)$ in the definition of $Z(t,C)$, we
obtain the $\zeta$-function of the curve $C/\mathbf{F}_q$, for distinct
cases of $i,j(\bmod~{2l})$, as below:

\begin{enumerate}
\renewcommand\labelenumi{\arabic{enumi}.}
\leftskip -.2pc
\item For $i \equiv j \equiv 0(\bmod~{2l})$,
\begin{align*}
\hskip -1.25pc \log Z(t,C) & = \sum_{n=1}^{\infty} \frac{(q^n + 1 -
(2l-1)(2l-2)(-1)^{ns} q^{n/2})t^n}{n} \\[.3pc]
\hskip -1.25pc &= \log \frac{1}{1-qt} + \log \frac{1}{1-t} -
(2l-1)(2l-2)\\[.3pc]
\hskip -1.25pc &\quad\ \times \log \frac{1}{1 - (-1)^sq^{1/2}t}.
\end{align*}
Hence
\begin{equation*}
\hskip -1.25pc Z(t, C) = \frac{(1 - (-1)^sq^{1/2}t)^{(2l-1)(2l-2)}}
{(1-t)(1-qt)}.
\end{equation*}

\item For $i,j \equiv 0(\bmod~{2}),\, i,j, i-j \not \equiv
0(\bmod~{2l})$,
\begin{alignat*}{2}
\hskip -1.25pc \log Z(t,C) & = & & \sum_{l | n} \frac{a_{2l}(n) t^n}{n}
+ \sum_{l \not\:|\, n} \frac{a_{2l}(n)t^n}{n} \\[.3pc]
\hskip -1.25pc & = & & \sum_{n=1}^{\infty} \frac{(q^{ln} + 1 -
(2l-1)(2l-2)(-1)^{lns} q^{ln/2})t^{ln}}{ln} \\[.3pc]
\hskip -1.25pc & & & + \sum_{n=1}^{\infty} \frac{(q^n + 1 - 2(-1)^{ns}
q^{n/2}) t^n}{n}\\[.3pc]
\hskip -1.25pc & & & - \sum_{n=1}^{\infty} \frac{(q^{ln} + 1 -
2(-1)^{lns} q^{ln/2}) t^{ln}}{ln} \\[.3pc]
\hskip -1.25pc & = & & \sum_{n=1}^{\infty} \frac{(q^n + 1 -
2(-1)^{ns}q^{n/2}) t^n}{n}\\[.3pc]
\hskip -1.25pc & & & - (4l-6) \sum_{n=1}^{\infty} \frac{(-1)^{lns}
q^{ln/2}t^{ln}}{n}\\[.3pc]
\hskip -1.25pc & = & & \log \frac{1}{1-qt} + \log \frac{1}{1-t} - 2 \log
\frac{1}{1 - (-1)^s q^{1/2}t}\\[.3pc]
\hskip -1.25pc & & & - (4l-6) \log \frac{1}{1 - (-1)^{ls}q^{l/2}t^l}.
\end{alignat*}
Hence
\begin{align*}
\hskip -1.25pc Z(t, C) &= \frac{(1 - (-1)^s q^{1/2}t)^2 (1 -
(-1)^{ls}q^{l/2}t^l)^{4l-6}}{(1-t)(1-qt)} \\[.3pc]
\hskip -1.25pc &= \frac{(1 - \theta t)^{4l-4} \prod_{r=1}^{l-1}(1 - \zeta^r \theta
t)^{4l-6}}{(1-t)(1-qt)}.
\end{align*}

\item For $i, j \equiv 0(\bmod~{2})$ and (i) $i \equiv 0,\, j \not
\equiv 0(\bmod~{2l})$, (ii) $i \not \equiv 0,\, j \equiv 0(\bmod~{2l})$,
(iii) $i,j \not\equiv 0,\, i \equiv j(\bmod~{2l})$,
\begin{alignat*}{2}
\hskip -1.25pc \log Z(t,C) &= & & \sum_{l | n} \frac{a_{2l}(n) t^n}{n}
+ \sum_{l \not\:|\, n} \frac{a_{2l}(n) t^n}{n} \\[.3pc]
\hskip -1.25pc &= & & \sum_{n=1}^{\infty} \frac{(q^{ln} + 1 -
(2l-1)(2l-2)(-1)^{lns}q^{ln/2}) t^{ln}}{ln} \\[.3pc]
\hskip -1.25pc & & & + \sum_{n=1}^{\infty} \frac{(q^n + 1 + 2(l-1) (-1)^{ns} q^{n/2})
t^n}{n} \displaybreak[0] \\[.3pc]
\hskip -1.25pc & & & - \sum_{n=1}^{\infty} \frac{(q^{ln} + 1 + 2(l-1)
(-1)^{lns}q^{ln/2}) t^{ln}}{ln} \\[.3pc]
\hskip -1.25pc & = & & \sum_{n=1}^{\infty} \frac{(q^n + 1 + 2(l-1)(-1)^{ns}
q^{n/2})t^n}{n} \\[.3pc]
\hskip -1.25pc & & & - (4l-4) \sum_{n=1}^{\infty} \frac{(-1)^{lns}q^{ln/2}t^{ln}}{n} \\[.3pc]
\hskip -1.25pc & = & & \log \frac{1}{1-qt} + \log \frac{1}{1-t} + 2(l-1) \log
\frac{1}{1 - (-1)^s q^{1/2}t} \\[.3pc]
\hskip -1.25pc & & & - (4l-4) \log \frac{1}{1 - (-1)^{ls}q^{l/2}t^l}.
\end{alignat*}
It follows that
\begin{align*}
\hskip -1.25pc Z(t, C) &= \frac{(1 - (-1)^{ls} q^{l/2}t^l)^{4l-4}} {(1-t)(1-qt)(1 -
(-1)^s q^{1/2}t)^{2l-2}} \\[.3pc]
\hskip -1.25pc &= \frac{(1 - \theta t)^{2l-2} \prod_{r=1}^{l-1} (1 - \zeta^r \theta
t)^{4l-4}}{(1-t)(1-qt)}.
\end{align*}

\item For (i) $i \equiv 0,\, j \equiv l(\bmod~{2l})$, (ii) $i \equiv
l,\, j \equiv 0(\bmod~{2l})$, and (iii) $i \equiv j \equiv
l(\bmod~{2l})$,
\begin{alignat*}{2}
\hskip -1.25pc \log Z(t,C) &= & & \sum_{2 | n} \frac{a_{2l}(n) t^n}{n} +
\sum_{2 \not\:|\, n} \frac{a_{2l}(n) t^n}{n} \\[.3pc]
\hskip -1.25pc &= & & \sum_{n=1}^{\infty} \frac{(q^{2n} + 1 -
(2l-1)(2l-2)(-1)^{2ns} q^{2n/2})t^{2n}}{2n} \\[.3pc]
\hskip -1.25pc & & & + \sum_{n=1}^{\infty} \frac{(q^n + 1 + 2(l-1)
(-1)^{ns} q^{n/2}) t^n}{n} \\[.3pc]
\hskip -1.25pc & & & - \sum_{n=1}^{\infty} \frac{(q^{2n} + 1 + 2(l-1)
(-1)^{2ns} q^{2n/2}) t^{2n}}{2n} \\[.3pc]
\hskip -1.25pc & = & & \sum_{n=1}^{\infty} \frac{(q^n + 1 +
2(l-1)(-1)^{ns} q^{n/2})t^n}{n} \\[.3pc]
\hskip -1.25pc & & & - 2l(l-1) \sum_{n=1}^{\infty} \frac{(-1)^{2ns}
q^{2n/2}t^{2n}}{n} \displaybreak[0] \\[.3pc]
\hskip -1.25pc & = & & \log \frac{1}{1-qt} + \log \frac{1}{1-t} +
2(l-1)\log \frac{1}{1 - (-1)^s q^{1/2} t} \\[.3pc]
\hskip -1.25pc & & & - 2l(l-1) \log \frac{1}{1 - (-1)^{2s} qt^2}.
\end{alignat*}
It follows that
\begin{align*}
\hskip -1.25pc Z(t, C) &= \frac{(1 - (-1)^{2s}qt^2)^{2l(l-1)}}
{(1-t)(1-qt)(1 - (-1)^s q^{1/2}t)^{2(l-1)}} \\[.4pc]
\hskip -1.25pc &= \frac{(1 - \theta t)^{2(l-1)^2} (1 + \theta
t)^{2l(l-1)}}{(1-t)(1-qt)}.
\end{align*}

\item For (i) $j \not\equiv 0(\bmod~{2}),\, i \equiv 0,\, j \not\equiv
l(\bmod~{2l})$, (ii) $i \not\equiv 0(\bmod~{2}),\, i \not \equiv l,\, j
\equiv 0(\bmod~{2l})$, and (iii) $i \not\equiv 0(\bmod~{2}),\, i \equiv
j,\, i \not\equiv l(\bmod~{2l})$,
\begin{alignat*}{2}
\hskip -1.25pc \log Z(t, C) &= & & \sum_{2l | n} \frac{a_{2l}(n) t^n}{n} + \sum_{2l
\not\:|\, n} \frac{a_{2l}(n) t^n}{n} \\[.3pc]
\hskip -1.25pc &= & & \sum_{n=1}^{\infty} \frac{(q^{2ln} + 1 -
(2l-1)(2l-2)(-1)^{2lns}q^{2ln/2}) t^{2ln}}{2ln} \\[.3pc]
\hskip -1.25pc & & & + \sum_{n=1}^{\infty} \frac{(q^n + 1 + 2(l-1) (-1)^{ns} q^{n/2})
t^n}{n} \\[.3pc]
\hskip -1.25pc & & & - \sum_{n=1}^{\infty} \frac{(q^{2ln} + 1 + 2(l-1) (-1)^{2lns}
q^{2ln/2}) t^{2ln}}{2ln} \\[.3pc]
\hskip -1.25pc & = & & \sum_{n=1}^{\infty} \frac{(q^n + 1 + 2(l-1)(-1)^{ns}
q^{n/2})t^n}{n} \\[.3pc]
\hskip -1.25pc & & & - 2(l-1) \sum_{n=1}^{\infty} \frac{(-1)^{2lns} q^{2ln/2}
t^{2ln}}{n} \\[.3pc]
\hskip -1.25pc & = & & \log \frac{1}{1-qt} + \log \frac{1}{1-t} + 2(l-1)\log \frac{1}{1
- (-1)^s q^{1/2} t} \\[.3pc]
\hskip -1.25pc & & & - 2(l-1)\log \frac{1}{1 - (-1)^{2ls}q^l t^{2l}}.
\end{alignat*}
This implies that
\begin{align*}
\hskip -1.25pc Z(t,C) & = \frac{(1 - (-1)^{2ls}q^l t^{2l})^{2l-2}} {(1-t)(1-qt)(1 -
(-1)^s q^{1/2}t)^{2l-2}} \displaybreak[0] \\[.4pc]
\hskip -1.25pc & = \frac{(1 + \theta^l t^l)^{2l-2} \prod_{r=1}^{l-1} (1 - \zeta^r
\theta t)^{2l-2}}{(1-t)(1-qt)}.
\end{align*}

\item For $i,j, i-j \not\equiv 0,l(\bmod~{2l})$ and (i) $i \not\equiv
j(\bmod~{2})$, (ii) $i,j \not\equiv 0(\bmod~{2})$,
\begin{alignat*}{2}
\hskip -1.25pc \log Z(t,C) & = & & \sum_{2l | n} \frac{a_{2l}(n) t^n}{n} + \sum_{2
\not\:|\, n,\; l | n} \frac{a_{2l}(n) t^n}{n} + \sum_{l \not\:|\, n}
\frac{a_{2l}(n) t^n}{n} \\[.3pc]
\hskip -1.25pc & = & & \sum_{n=1}^{\infty} \frac{(q^{2ln} + 1 - (2l-1)(2l-2)(-1)^{2lns}
q^{2ln/2}) t^{2ln}}{2ln} 
\end{alignat*}
\begin{alignat*}{2}
\hskip -1.25pc & & & + \sum_{n=1}^{\infty} \frac{(q^{ln} + 1 + 2(l-1) (-1)^{lns}
q^{ln/2}) t^{ln}}{ln} \\[.3pc]
\hskip -1.25pc & & & - \sum_{n=1}^{\infty} \frac{(q^{2ln} + 1 + 2(l-1) (-1)^{2lns}
q^{2ln/2}) t^{2ln}}{2ln} \\[.3pc]
\hskip -1.25pc & & & + \sum_{n=1}^{\infty} \frac{(q^n + 1 - 2(-1)^{ns} q^{n/2})
t^n}{n}\\[.3pc]
\hskip -1.25pc & & & - \sum_{n=1}^{\infty} \frac{(q^{ln} + 1 - 2(-1)^{lns} q^{ln/2})
t^{ln}}{ln} \\[.3pc]
\hskip -1.25pc & = & & \sum_{n=1}^{\infty} \frac{(q^n + 1 -2(-1)^{ns} q^{n/2})
t^n}{n}\\[.3pc]
\hskip -1.25pc & & & - (2l-2) \sum_{n=1}^{\infty} \frac{(-1)^{2lns}
q^{2ln/2}t^{2ln}}{n} + 2 \sum_{n=1}^{\infty} \frac{(-1)^{lns}q^{ln/2}
t^{ln}}{n}\\[.3pc]
\hskip -1.25pc & = & & \log \frac{1}{1-qt} + \log \frac{1}{1-t} -2\log \frac{1}{1 -
(-1)^s q^{1/2} t} \\[.3pc]
\hskip -1.25pc & & & - (2l-2)\log \frac{1}{1 - (-1)^{2ls}q^l t^{2l}} + 2 \log\frac{1}{1
- (-1)^{ls}q^{l/2} t^l}.
\end{alignat*}
Thus we obtain
\begin{align*}
\hskip -1.25pc Z(t, C) &= \frac{(1 - (-1)^s q^{1/2}t)^2 (1 -
(-1)^{2ls}q^l t^{2l})^{2l-2}}{(1-t)(1-qt)(1 - (-1)^{ls}q^{l/2}t^l)^2} \\[.4pc]
\hskip -1.25pc &= \frac{(1 - \theta t)^2 (1 - \theta^l t^l)^{2l-4} (1 +
\theta^l t^l)^{2l-2}}{(1-t)(1-qt)}.
\end{align*}

\item For (i) $i \equiv l,\, j \not\equiv 0,l(\bmod~{2l})$, (ii) $i
\not\equiv 0,l,\, j \equiv l(\bmod~{2l})$, and (iii) $i,j \not\equiv
0,l,\, i-j \equiv l(\bmod~{2l})$,
\begin{alignat*}{2}
\hskip -1.25pc \log Z(t,C) &= & & \sum_{2l | n} \frac{a_{2l}(n) t^n}{n} + \sum_{2
\not\:|\, n,\; l | n}\frac{a_{2l}(n)t^n}{n}\\[.3pc]
& & & + \sum_{2 | n,\; l \not\:|\,
n}\frac{a_{2l}(n)t^n}{n} + \sum_{2 \not\:|\, n,\; l \not\:|\,
n}\frac{a_{2l}(n)t^n}{n} \\[.3pc]
\hskip -1.25pc & = & & \sum_{n=1}^{\infty} \frac{(q^{2ln} + 1 - (2l-1)(2l-2)
(-1)^{2lns}q^{2ln/2})t^{2ln}}{2ln} \\[.3pc]
\hskip -1.25pc & & & + \sum_{n=1}^{\infty} \frac{(q^{ln} + 1 + 2(l-1)(-1)^{lns}
q^{ln/2})t^{ln}}{ln} \\[.3pc]
\hskip -1.25pc & & & - \sum_{n=1}^{\infty} \frac{(q^{2ln} + 1 + 2(l-1) (-1)^{2lns}
q^{2ln/2})t^{2ln}}{2ln} 
\end{alignat*}
\begin{alignat*}{2}
\hskip -1.25pc & & & + \sum_{n=1}^{\infty} \frac{(q^{2n} + 1 + 2(l-1)(-1)^{2ns}
q^{2n/2}) t^{2n}}{2n} \\[.3pc]
\hskip -1.25pc & & & - \sum_{n=1}^{\infty} \frac{(q^{2ln} + 1 + 2(l-1) (-1)^{2lns}
q^{2ln/2}) t^{2ln}}{2ln} \\[.3pc]
\hskip -1.25pc & & & + \sum_{n=1}^{\infty} \frac{(q^n + 1 - 2(-1)^{ns} q^{n/2}) t^n}{n}
\displaybreak[0] \\[.3pc]
\hskip -1.25pc & & & - \sum_{n=1}^{\infty} \frac{(q^{2n} + 1 - 2(-1)^{2ns} q^{2n/2})
t^{2n}}{2n} \\[.3pc]
\hskip -1.25pc & & & - \sum_{n=1}^{\infty} \frac{(q^{ln} + 1 - 2(-1)^{lns} q^{ln/2})
t^{ln}}{ln} \\[.3pc]
\hskip -1.25pc & & & + \sum_{n=1}^{\infty} \frac{(q^{2ln} + 1 - 2(-1)^{2lns} q^{2ln/2})
t^{2ln}}{2ln} \\[.3pc]
\hskip -1.25pc & = & & \sum_{n=1}^{\infty} \frac{(q^n + 1 - 2(-1)^{ns} q^{n/2}) t^n}{n}
+ 2 \sum_{n=1}^{\infty} \frac{(-1)^{lns} q^{ln/2} t^{ln}}{n} \\[.3pc]
\hskip -1.25pc & & & - (2l-1) \sum_{n=1}^{\infty} \frac{(-1)^{2lns} q^{2ln/2}
t^{2ln}}{n} + l\sum_{n=1}^{\infty} \frac{(-1)^{2ns} q^{2n/2}t^{2n}}{n}\\[.3pc]
\hskip -1.25pc & = & & \log \frac{1}{1-qt} + \log \frac{1}{1-t} - 2\log \frac{1}{1 -
(-1)^s q^{1/2}t} \displaybreak[0]\\[.3pc]
\hskip -1.25pc & & & + 2\log \frac{1}{1 - (-1)^{ls} q^{l/2} t^l} - (2l-1)\log
\frac{1}{1 - (-1)^{2ls}q^l t^{2l}} \\[.3pc]
\hskip -1.25pc & & & + l\log \frac{1}{1 - (-1)^{2s} qt^2}.
\end{alignat*}
Hence we obtain
\begin{align*}
\hskip -1.25pc Z(t,C) & = \frac{(1 - (-1)^s q^{1/2}t)^2 (1 - (-1)^{2ls}q^l
t^{2l})^{2l-1}}{(1-t)(1-qt)(1 - (-1)^{ls}q^{l/2}t^l)^2 (1 -
(-1)^{2s}qt^2)^{l}} \\[.5pc]
\hskip -1.25pc & = \frac{(1 - \theta^l t^l)^{2l-3} (1 + \theta^l
t^l)^{2l-1}}{(1-t)(1-qt) (1 - \theta t)^{l-2} (1 + \theta t)^{l}} \\[.5pc]
\hskip -1.25pc & = \frac{((1 - \theta t)(1 + \theta t))^{l-1} \prod_{r=1}^{l-1} ((1 -
\zeta^r \theta t)^{2l-3}(1 + \zeta^r \theta t)^{2l-1})} {(1-t)(1-qt)}.
\end{align*}
Hence the theorem.\hfill $\Box$\vspace{-1.5pc}
\end{enumerate}
\end{proof}

The curve $C$ in Theorem~\ref{thm-1} is non-singular of degree $2l$;
hence it has genus $g = (2l-1)(2l-2)/2$. From the expressions for $P(t)$
in Theorem~\ref{thm-1}, it is clear that in each case for
$i,j(\bmod~{2l})$,
\begin{itemize}
\leftskip -1pc
\item $P(t)$ is a polynomial of degree $2g = (2l-1)(2l-2)$ of the form
$P(t) = \prod_{k=1}^{2g} (1 - \alpha_k t)$, and the $\alpha_k$'s are
algebraic integers equal to $\pm q^{1/2} \zeta^r$, $0 \leq r \leq l-1$.
Thus $| \alpha_k | = q^{1/2}$ for $1 \leq k \leq 2g$.

\item We may pair the $\alpha_k$'s in such a way that $\alpha_k
\alpha_{g+k} = q$ for $1 \leq k \leq g$. The polynomial $P(t)$ has
integer coefficients (as $q^{1/2}$ is an integer), constant term 1, and
leading term $q^g$ (since $\prod_{k=1}^{2g} \alpha_k = q^g$ by the above
pairing).
\end{itemize}
This corroborates the Weil conjectures for the concrete class of curves
$C/\mathbf{F}_q$ considered in Theorem~\ref{thm-1}.

Given a non-singular projective curve $X$ defined over a finite field
$k$, it is well-known that the class number $h$ of the function field of
$X/k$ satisfies the relation $h=P(1)$, where $P(t)$ is the polynomial
that appears in the numerator of the $\zeta$-function of the\break curve
$X/k$.

For the curve $C/\mathbf{F}_q$ considered in Theorem~\ref{thm-1}, and
the $\zeta$-function(s) obtained therein, we may thus substitute $t=1$
in the expressions for the polynomials $P(t)$, to obtain the class
number(s) of the associated function field(s) in Theorem~\ref{thm-2}
below.

\begin{theor}[\!]\label{thm-2}
Consider the projective curve $C \colon aY^{2l} = bX^{2l} + cZ^{2l}$
$(abc \neq 0)$ defined over the finite field $\mathbf{F}_q${\rm ,} with
notations as in Theorem~{\rm \ref{thm-1}}. Set $q_0 = p^f$ and $u =
\sqrt{q_0}$. Thus $u$ is an integer{\rm ,} $q = q_0^s${\rm ,} and $\sqrt{q} = u^s$.
For each $s \geq 1${\rm ,} let $K_s$ denote the function field of the curve
$C/\mathbf{F}_q${\rm ,} $q=p^{fs}${\rm ,} and let $h_s$ denote its class number. Let
$h_s = h_1$ for $s$ odd{\rm ,} and $h_s = h_2$ for $s$ even. Substituting $h_s
= P(1)${\rm ,} we obtain the class numbers $h_1$ and $h_2${\rm ,} for the seven
distinct cases in Theorem~{\rm \ref{thm-1},} as below{\rm :}
\begin{enumerate}
\renewcommand\labelenumi{\rm \arabic{enumi}.}
\leftskip -.2pc
\item For $i \equiv j \equiv 0(\bmod~{2l})${\rm ,}
\begin{equation*}
\hskip -1.25pc h_1 = (u^s + 1)^{(2l-1)(2l-2)}, \quad h_2 = (u^s - 1)^{(2l-1)(2l-2)}.
\end{equation*}
\item For $i,j \equiv 0(\bmod~{2}), i,j, i-j \not\equiv
0(\bmod~{2l})${\rm ,}
\begin{equation*}
\hskip -1.25pc h_1 = (u^s + 1)^2 (u^{ls} + 1)^{4l-6}, \quad h_2 = (u^s - 1)^2 (u^{ls}
- 1)^{4l-6}.
\end{equation*}
\item For $i, j \equiv 0(\bmod~{2})$ and {\rm (i)} $i \equiv 0,\, j
\not\equiv 0(\bmod~{2l})${\rm ,} {\rm (ii)} $i \not\equiv 0,\, j \equiv
0(\bmod~{2l})${\rm ,} {\rm (iii)}~$i,j \not\equiv 0,\, i \equiv
j(\bmod~{2l})${\rm ,}
\begin{equation*}
\hskip -1.25pc h_1 = (u^{ls} + 1)^{4l-4}/(u^s + 1)^{2l-2}, \quad h_2 = (u^{ls} -
1)^{4l-4}/(u^s - 1)^{2l-2}.
\end{equation*}

$\left.\right.$\vspace{-1pc}

\item For {\rm (i)} $i \equiv 0,\, j \equiv l(\bmod~{2l})${\rm ,} {\rm
(ii)} $i \equiv l,\, j \equiv 0(\bmod~{2l})${\rm ,} and {\rm (iii)} $i
\equiv j \equiv l(\bmod~{2l})${\rm ,}
\begin{equation*}
\hskip -1.25pc h_1 = (u^s + 1)^{2(l-1)^2} (u^s - 1)^{2l(l-1)}, \quad h_2
= (u^s - 1)^{2(l-1)^2} (u^s + 1)^{2l(l-1)}.
\end{equation*}

$\left.\right.$\vspace{-1.5pc}

\item For {\rm (i)} $j \not\equiv 0(\bmod~{2}),\, i \equiv 0,\, j
\not\equiv l(\bmod~{2l})${\rm ,} {\rm (ii)} $i \not\equiv
0(\bmod~{2}),\, i \not\equiv l,\, j \equiv 0(\bmod~{2l})${\rm ,} and
{\rm (iii)} $i \not\equiv 0(\bmod~{2}),\, i \equiv j,\, i \not \equiv
l(\bmod~{2l})${\rm ,}
\begin{equation*}
\hskip -1.25pc h_1 = ((u^{2ls} - 1)/(u^s + 1))^{2l-2}, \quad h_2 =
((u^{2ls} -1)/(u^s - 1))^{2l-2}.
\end{equation*}
\item For $i,j, i-j \not\equiv 0,l(\bmod~{2l})$ and {\rm (i)} $i
\not\equiv j(\bmod~{2})${\rm ,} {\rm (ii)} $i,j \not\equiv
0(\bmod~{2})${\rm ,}
\begin{align*}
\hskip -1.25pc h_1 &= (u^s + 1)^2 (u^{ls} + 1)^{2l-4} (u^{ls} -1)^{2l-2},\\[.2pc]
\hskip -1.25pc h_2 &= (u^s -1)^2 (u^{ls} - 1)^{2l-4} (u^{ls} + 1)^{2l-2}.
\end{align*}
\item For {\rm (i)} $i \equiv l,\, j \not\equiv 0,l(\bmod~{2l})${\rm ,}
{\rm (ii)} $i \not\equiv 0,l,\, j \equiv l(\bmod~{2l})${\rm ,} and
{\rm (iii)} $i,j \not\equiv 0,l,\, i-j \equiv l(\bmod~{2l})${\rm ,}
\begin{align*}
\hskip -1.25pc h_1 &= (u^{ls} + 1)^{2l-3} (u^{ls} - 1)^{2l-1}/ (u^s +
1)^{l-2} (u^s - 1)^{l},\\[.2pc]
\hskip -1.25pc h_2 &= (u^{ls} - 1)^{2l-3} (u^{ls} + 1)^{2l-1}/ (u^s -
1)^{l-2} (u^s + 1)^{l}.
\end{align*}
\end{enumerate}
\end{theor}

For the class of curves $C/\mathbf{F}_q$ considered in
Lemma~\ref{lem-1}, the explicit form of the polynomials $P(t)$
in the $\zeta$-function(s) for $C/\mathbf{F}_q$ were obtained in
(\cite{an-1}, Theorem~8).  Substituting $t=1$ in these expressions
for $P(t)$, we obtain the class number(s) of the associated
function field(s) in Theorem~\ref{thm-3} below.

\begin{theor}[\!]\label{thm-3}
Consider the projective curve $C \colon aY^l = bX^l + cZ^l$ $(abc \neq
0)$ defined over the finite field $\mathbf{F}_q${\rm ,} with notations
as in Lemma~{\rm \ref{lem-1}}. Let $q_0 = p^f$ and $u = \sqrt{q_0}$.
Thus $u$ is an integer{\rm ,} $q = q_0^s${\rm ,} and $\sqrt{q} = u^s$.
Let $\theta = (-1)^s q^{1/2}$ and let $\zeta$ be a primitive {\rm
(}complex{\rm )} $l$-th root of unity. For each $s \geq 1${\rm ,} let
$K_s$ denote the function field of the curve $C/\mathbf{F}_q${\rm ,}
$q=p^{fs}${\rm ,} and let $h_s$ denote its class number. Let $h_s = h_1$
for $s$ odd{\rm ,} and $h_s = h_2$ for $s$ even. Substituting $h_s =
P(1)$ in the expressions for the polynomial $P(t)$ in the
$\zeta$-function{\rm (}s{\rm )} of the curve $C/\mathbf{F}_q$ {\rm (}cf.
{\rm (\cite{an-1},} Theorem~$8${\rm ),} reproduced below{\rm ),} we
obtain the class numbers $h_1$ and $h_2${\rm ,} for the distinct cases
in Lemma~{\rm \ref{lem-1},} as below{\rm :}
\begin{enumerate}
\renewcommand\labelenumi{\rm \arabic{enumi}.}
\leftskip -.2pc
\item For $i,j \equiv 0(\bmod~{l})${\rm ,}
\begin{align*}
\hskip -1.25pc &P(t) = (1 - \theta t)^{(l-1)(l-2)},\\[.2pc]
\hskip -1.25pc &h_1 = (u^s + 1)^{(l-1)(l-2)}, \quad h_2 = (u^s - 1)^{(l-1)(l-2)}.
\end{align*}
\item For {\rm (i)} $i \equiv 0(\bmod~{l}), j \not\equiv
0(\bmod~{l})${\rm ,} {\rm (ii)} $i \not\equiv 0(\bmod~{l}), j \equiv
0(\bmod~{l})${\rm ,} and {\rm (iii)}~$i,j \not\equiv 0(\bmod~{l}), i
\equiv j(\bmod~{l})${\rm ,}
\begin{align*}
\hskip -1.25pc &P(t) = \prod_{r=1}^{l-1} (1 - \zeta^r \theta t)^{l-2},\\[.2pc]
\hskip -1.25pc &h_1 = ((u^{ls} + 1)/(u^s + 1))^{l-2}, \quad h_2 =
((u^{ls} - 1)/(u^s - 1))^{l-2}.
\end{align*}
\item For $i, j, i-j \not\equiv 0(\bmod~{l})${\rm ,}
\begin{align*}
\hskip -1.25pc &P(t) = (1 - \theta t)^{l-1} \prod_{r=1}^{l-1} (1 - \zeta^r \theta
t)^{l-3},\\[.3pc]
\hskip -1.25pc &h_1 = (u^s + 1)^2 (u^{ls} + 1)^{l-3}, \quad h_2 = (u^s - 1)^2 (u^{ls}
- 1)^{l-3}.
\end{align*}
\end{enumerate}
\end{theor}

Consider now the projective curves $C/\mathbf{F}_q$ in
Theorems~\ref{thm-2} and \ref{thm-3} as defined over some fixed base
field $\mathbf{F}_q$, $q=p^{fs_0}$, with associated function field
$K_{s_0}$. Then for each $s \geq 1$, the function fields $K_{ss_0}$ of
the curves $C/\mathbf{F}_{q^s}$, are isomorphic to the constant field
extensions $K_{s_0} \cdot \mathbf{F}_{q^s}$ of the function field $K_{s_0}$.
In this case, the results in Theorems~\ref{thm-2} and~\ref{thm-3}
provide concrete information on the growth of class numbers $h_{ss_0}$
($s \geq 1$) for the constant field extensions
$K_{s_0} \cdot \mathbf{F}_{q^s}/K_{s_0}$, for each class of curves. Note that
in this consideration, two cases arise: (i) for $s_0$ odd, the results
for both $h_1$ and $h_2$ come into picture, while (ii) for $s_0$ even,
only the results for $h_2$ are required.

\begin{conremar}$\left.\right.$
\begin{enumerate}
\renewcommand\labelenumi{\rm \arabic{enumi}.}
\leftskip -.2pc
\item In Theorem~\ref{thm-1}, for each distinct case, the roots of the
polynomial $P(t)$, all of which lie on the circle $| z | = q^{-1/2}$ in
the complex plane, are {\it not} uniformly distributed on this
circle. In each case, the roots are of the form $\beta = \xi q^{-1/2}$,
where $\xi$ is a complex $2l$-th root of unity. This is similarly the
case for the class of curves $aY^l = bX^l + cZ^l$ considered in
(\cite{an-1}, Theorem~8).

\item For each distinct case in Theorems~\ref{thm-2} and
\ref{thm-3}, the class number is a polynomial in $\sqrt{q}$,
of degree $2g$, with integral coefficients and constant term 1, where
$g$ is the genus of the curve $C/\mathbf{F}_q$.  (The genus of the
curve $C/\mathbf{F}_q$ in Theorem~\ref{thm-3} is $g = (l-1)(l-2)/2$.)

\item The polynomial $P(t)$ in the $\zeta$-function of a maximal curve
$C/\mathbf{F}_q$ is of the form
\begin{equation*}
\hskip -1.25pc P(t) = (1 + q^{1/2}t)^{2g},
\end{equation*}
and that of a minimal curve $C/\mathbf{F}_q$ is of the form
\begin{equation*}
\hskip -1.25pc P(t) = (1 - q^{1/2}t)^{2g}.
\end{equation*}
These expressions follow easily from the Weil conjectures applied to the
expression for $a(1, C)$ in \S\ref{sec2}. Conversely, given a
non-singular projective curve $C/\mathbf{F}_q$, such that the polynomial
$P(t)$ in its $\zeta$-function has the above form(s), one sees that the
curve $C$ is maximal (resp. minimal) over $\mathbf{F}_q$.

\hskip 1pc From the explicit expressions for the polynomial $P(t)$ in the
$\zeta$-function(s) of the projective curve $aY^e = bX^e + cZ^e$ defined
over $\mathbf{F}_q$ (cf. \cite{an-1}, Theorem~8) (for $e=l$) and
Theorem~\ref{thm-1} (for $e=2l$)), it is clear that these curves are
maximal (or minimal) over $\mathbf{F}_q$ precisely when the coefficients
$a,b,c$ reduce to the case $a=b=c=1$ corresponding to the Fermat curves.

\item In (\cite{katre}, Proposition 2), the author has stated results
(to appear) for the polynomials $P(t)$ in the $\zeta$-functions of the
projective curves $aY^e = bX^e + cZ^e$ ($abc \neq 0$) defined over
finite fields $\mathbf{F}_q$, $q = p^\alpha \equiv 1(\bmod~{e})$, for
integers $e \geq 3$ and primes $p$ such that ${\rm ord}\ p(\bmod~{e})$
is even. These results generalize the results obtained for the
polynomials $P(t)$ in (\cite{an-1}, Theorem~8) and Theorem~\ref{thm-1}
of this paper.\vspace{-.5pc}
\end{enumerate}
\end{conremar}

\section*{Acknowledgements}

The author thanks CSIR, India for support received in the form
of a research fellowship, when this paper was conceived.

\end{document}